# On Maps with a Single Zigzag

Sóstenes Lins [a], Valdenberg Silva [b]

[a] *Dept. Matemática da UFPE - Recife - Brazil*
[b] *Dept. Matemática da UFSE - Aracaju - Brazil*

**Abstract**

If a graph $G_M$ is embedded into a closed surface $S$ such that $S \backslash G_M$ is a collection of disjoint open discs, then $M = (G_M, S)$ is called a *map*. A *zigzag* in a map $M$ is a closed path which alternates choosing, at each star of a vertex, the leftmost and the rightmost possibilities for its next edge. If a map has a single zigzag we show that the cyclic ordering of the edges along it induces linear transformations, $c_P$ and $c_{P\sim}$ whose images and kernels are respectively the cycle and bond spaces (over $GF(2)$) of $G_M$ and $G_D$, where $D = (G_D, S)$ is the dual map of $M$. We prove that $Im(c_P \circ c_{P\sim})$ is the intersection of the cycle spaces of $G_M$ and $G_D$, and that the dimension of this subspace is connectivity of $S$. Finally, if $M$ has also a single face, this face induces a linear transformation $c_D$ which is invertible: we show that $c_D^{-1} = c_{P\sim}$.

## 1 Introduction: Combinatorial Maps

A *topological map* $M^t = (G, S)$ is an embedding of a graph $G$ into a closed surface $S$ such that $S \backslash G$ is a collection of disjoint open disks, called *faces*. By going around the boundary of a face and recalling the edges traversed we define a *facial path* of $M^t$, which is a closed path in $G$. Note that a facial path is obtained starting in an edge and by choosing at each vertex always the rightmost or always the leftmost possibility for the next edge. If we alternate the choice, then the result is a *zigzag path*, or simply a *zigzag*. Even if the surface is non-orientable these left-right choices are well defined, because they are local. For more background on graphs embedded into surfaces see [Giblin, 1977]. To make our objects less dependent of topology we use a combinatorial counterpart for topological maps introduced in [Lins, 1982]. A *combinatorial map* or simply a *map* $M$ is an ordered triple $(C_M, v_M, f_M)$ where: $(i)$ $C_M$ is a connected finite cubic graph; $(ii)$ $v_M$ and $f_M$ are disjoint perfect matchings in $C_M$, such that each component of the subgraph of $C_M$ induced by $v_M \cup f_M$ is a *polygon* (i.e. a non-empty connected subgraph with all the vertices having two incident edges) with 4 edges and it is called an $M$-*square*.

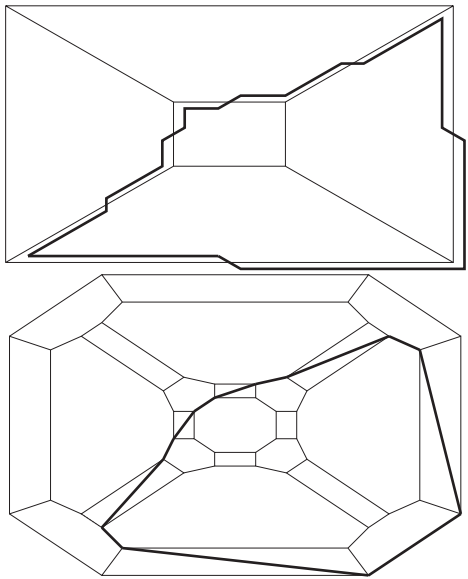

Fig. 1: A zigzag and its corresponding $z$-gon

From the above definition, it follows that $C_M$ may contain double edges but not loops. A third perfect matching in $C_M$ is $E(C_M) - (v_M \cup f_M)$ and is denoted by $a_M$. The set of diagonals of the $M$-squares, denoted by $z_M$, is a perfect matching in the complement of $C_M$. The edges in $v_M, f_M, z_M, a_M$ are called respectively $v_M$-edges, $f_M$-edges, $z_M$-edges, $a_M$-edges. The graph $C_M \cup z_M$ is denoted by $Q_M$, and is a regular graph of valence 4. A component induced by $a_M \cup v_M$ is a polygon with an even number of vertices and it is called a $v$-gon. Similarly, we define an $f$-gon, and a $z$-gon, by replacing $v$ for $f$ and $v$ for $z$. Clearly, the $f$-gons and $z$-gons of $C_M$ correspond to the facial paths and the zigzags of $M^t$. To avoid the use of colors the $M$-squares are presented in the pictures as rectangles in which the short sides $(s)$ are $v_M$-edges, the long sides $(\ell)$ are $f_M$-edges and the diagonals $(d)$ are $z_M$ edges. An $M$-*rectangle with diagonals* or simply an $M$-*rectangle* (being understood that the diagonals are present) is a component induced by $v_M, f_M, z_M$. The set of $M$-rectangles is denoted by $R$. If $\pi$ is a permutation of the symbols $s\ell d$, and $R' \subseteq R$ subset of rectangles of $M$, then $M(R' : \pi)$ denotes the map obtained from $M$ by permuting the short sides, the long sides and the diagonals according to $\pi$ in all $r \in R'$. Let $M(r : \pi)$ denote $M(\{r\} : \pi)$. The *dual map* of $M$ is the map $D = M(R : \ell s d)$; $D$ and $M$ have the same $z$-gons and the $v$-gons and $f$-gons interchanged. The *phial map* of $M$ is the map $P = M(R : d\ell s)$; $P$ and $M$ have the same $f$-gons and the $v$-gons and $z$-gons interchanged. The *antimap* of $M$ is the map $M\textasciitilde = M(R : sd\ell)$; $M$ and $M\textasciitilde$ have the same $v$-gons and the $f$-gons and $z$-gons interchanged. The pairs $(M, D), (M, P), (M, M\textasciitilde)$ constitute the *map dualities* introduced in [Lins, 1982]. The dual of $P$ is $D\textasciitilde$ and the dual of $M^\sim$ is $P^\sim$. Let $\Omega(M) = \{M, D, P, M^\sim, D^\sim, P^\sim\}$ and $\Omega^\star(M) = \{M(R' : \pi) \mid R' \subseteq R, \pi$ permutation of $s\ell d\}$. Note that $\Omega(M) \subseteq \Omega^\star(M)$ and that any member of $\Omega^\star(M)$ has $R$ as its set of rectangles.

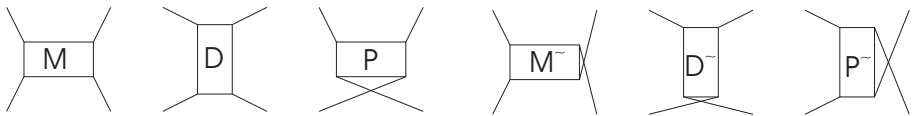

Fig. 2: How a neighborhood of each rectangle is modified in the members of $\Omega(M)$

Given a map $M$ and its dual $D$, there exists a closed surface, denoted by $\text{Surf}(M, D)$ where $C_M = C_D$ naturally embeds. Consider the $v$-gons, the $f$-gons and the $M$-squares bounding disjoint closed disks. Each edge of $C_M$ occurs twice in the boundary of this collection of disks. Identify the collection of disks along the two occurrences of each edge. The result is a closed surface and $C_M$ is *faithfully embedded on it*, meaning that the boundaries of the faces are *bicolored poly*gons or *bigons*. Similarly, there are surfaces $\text{Surf}(D^\sim, P)$ and $\text{Surf}(P^\sim, M^\sim)$.

We define a function $\psi$ which turns out to be a bijection from the set of maps onto the set of $t$-maps. We denote $\psi(M)$ by $M^t$. Given a map $M$, to obtain $M^t$ we proceed as follows. Consider the $t$-map $(C_M, S)$, where $S = \text{Surf}(M, D)$, given by the faithful embedding of $M$. The $v$-gons, the $f$-gons and the $M$-squares are boundaries of (closed, in this case) disks embedded (and forming) the surface $S(M)$. Shrink to a point the disjoint closed disks bounded by $v$-gons. The $M$-squares, then, become bounding digons. Shrink each such bounding digon to a line, maintaining unaffected its vertices. With these contractions, effected in $S$, $t$-map $(C_M, S)$ becomes, by definition, $M^t = (G_M, S)$. Graph $G_M$ is called the *graph induced* by $M$. A combinatorial description of $G_M$ can be given as follows: the vertices of $G_M$ are the $v$-gons of $M$; its edges are the squares of $M$; the two ends of an edge of $G_M$ are the two $v$-gons (which may coincide and the edge is a *loop*) that contain the $v_M$-edges of the corresponding $M$-square. It is evident that $\psi$ is inversible: given a $t$-map we replace each edge by a bounding digon in its surface, and then expand each vertex to a disc in order to obtain a cellular embedding of a cubic graph. Therefore, $\psi^{-1}$ is well-defined; in fact, it is the dual of a useful construction in topology, namely, barycentric division. Thus, $\psi$ is a bijection from the set of maps onto the set of $t$-maps. It can be observed that $\psi$ induces a bijection from the set of $M$-rectangles onto the set of edges of $G_M$. We use this bijection to identify the sets $R$ and $E(G_M)$. Via $R$, which is invariant for the members of $\Omega^\star(M)$, we identify $E(G_M)$ and $E(G_{M'})$ for $M' \in \Omega^\star(M)$. Denote these identified sets of edges by $E$.

## 2  Absorption Property on Maps

The members of $\Omega(M)$ induce three distinct graphs: $G_M = G_{M\sim}$, $G_D = G_{D\sim}$, $G_P = G_{P\sim}$. We give to the power set $\mathcal{E} = \{E' \mid E' \subseteq E\}$ a vector space structure over the field $GF(2)$ by defining the sum of subsets of edges $\sum\{A_i \mid 1 \leq i \leq n\}$ to mean the subset of $\bigcup\{A_i \mid 1 \leq i \leq n\}$ formed by the elements which occurs in an odd number of $A_i$'s. The *bond space* and the *cycle space* [Bondy and Murty, 1976], [Godsil and Royle, 2001] of $G_M$ are denoted respectively by $\mathcal{V}$ and $\mathcal{V}^\perp$. For connected $C_M$, they are vector subspaces of $\mathcal{E}$ whose dimensions are $v-1$ and $|E|-v+1$, where $v$ is the number of vertices of $G_M$, or number of $v$-gons of $M$. Similarly, let $\mathcal{F}$ and $\mathcal{F}^\perp$ denote the bond and cycle space of $G_D$ and $\mathcal{Z}$ and $\mathcal{Z}^\perp$ the bond and cycle space of $G_P$. The dimensions of $\mathcal{F}$ and $\mathcal{F}^\perp$ are $f-1$ and $|E|-f+1$, where $f$ is the number of $f$-gons of $M$, or *faces* of $M^t = (G_M, \text{Surf}(M, D))$. The dimensions of $\mathcal{Z}$ and $\mathcal{Z}^\perp$ are $z-1$ and $|E|-z+1$, where $z$ is the number of $z$-gons of $M$. Each $z$-gon in $M$ corresponds to a *zigzag* in $M^t = (G_M, \text{Surf}(M, D))$. If $G$ is a graph and $W \subseteq \mathcal{V}(G)$ then $\delta_G(W) = \{x \in E(G) \mid x \text{ has an end in } W \text{ and an end in } V(G)\setminus W\}$. The subset of edges $\delta_G(W)$ is called the *bond of $W$* in graph $G$.

**Theorem 1 (Absorption Property)** *For an arbitrary map $M$ we have*

$$(a)\ \mathcal{V} \cap \mathcal{F} \subseteq \mathcal{Z},\ (b)\ \mathcal{F} \cap \mathcal{Z} \subseteq \mathcal{V},\ (c)\ \mathcal{Z} \cap \mathcal{V} \subseteq \mathcal{F}.$$

**Proof:** It is enough to prove the result in case $(a)$. For $(b)$ and $(c)$ we would use maps $D^\sim$ and $P^\sim$, respectively. Take an element $X \in \mathcal{V} \cap \mathcal{F}$. It follows that there exist $V' \subseteq V(G_M)$ and $F' \subseteq V(G_D)$ such that $X = \delta_{G_M}(V')$ and $X = \delta_{G_D}(F')$. Denote by $V''$ the subgraph of $C_M$ consisting of the disjoint $v$-gons corresponding to the vertices in $V'$. Denote by $F''$ the subgraph of $C_M$ consisting of the disjoint $f$-gons corresponding to the vertices in $F'$. Let $Z'''$ the subgraph induced by the symmetric difference of $V''$ and $F''$. Since $C_M$ is cubic, $Z'''$ is a collection of disjoint polygons. The crucial observation is that any rectangle has precisely one opposite pair of vertices in $F'''$ if and only if this rectangle has one $v_M$-edge in $V''$ and one $f_M$-edge in $F''$, namely, if and only if the rectangle corresponds to an edge in $X$. Therefore, each component of $Z'''$ is a polygon with $3 \times k$ edges which can be factored by $k$ subpaths of length 3 having an $a_M$-edge, a $v_M$-edge and an $f_M$-edge (the last two not necessarily in this order). Let $Z''$ obtained from $Z'''$ by replacing each adjacent pair of $v_M$-, $f_M$- edges by the $z_M$-edge forming a triangle with them in $Q_M$. Clearly, $Z''$ is a collection of disjoint $z$-gons corresponding to a $Z' \subseteq V(G_P)$ satisfying $\delta_{G_M}(Z') = X$  ∎

# 3 Maps with a single zigzag

Maps with a single zigzag are related to the *Gauss code problem*. For this problem in the plane see [Shank, 1975], [Lovasz and Max, 1976], [Rosenstiehl, 1976]. The present paper has its motivation in trying to generalize the result of [Lins,, Richter, and Shank, 1987] on the 2-face colorable Gauss code problem in the projective plane to surfaces of higher genus. In fact, by using the algebraic theory we are about to develop we can solve the 2-face colorable Gauss code problem in the Klein bottle. The paper presenting this solution is under preparation. Thus, maps with a single zigzag induce useful algebraic structures which are the main theme of this work. Here is an example of a $t$-map having a single zigzag.

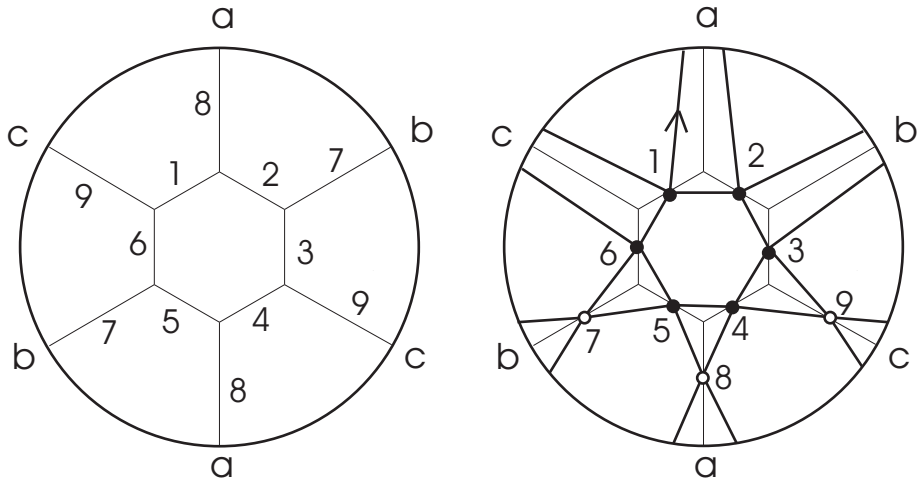

Fig. 3: An embedding of $K_{3,3}$ in the projective plane having a single zigzag

The cyclic sequence of edges visited in the single zigzag is

$$P = (1, 8, 5, 6, 9, 4, 5, 7, 3, 4, -8, 2, 3, -9, 1, 2, -7, 6).$$

This can be easily followed in the *medial map* [Godsil and Royle, 2001] on the right. The edges become small circles and, to obtain the faithful embedding of $C_M$, it is enough to deform each such circle to a rectangle. The direction of the first occurrence of an edge defines its orientation. Edges 1,2,3,4,5,6 are traversed twice in the positive direction (they correspond to black circles in the medial map) and edges 7,8,9 are traversed once in the positive direction and once in the negative direction (they correspond to white circles). The reason for the notation $P$ is that the signed cyclic sequence defines the phial map $P$ (whence all maps in $\Omega(M)$) and vice-versa, the phial defines the sequence.

Given a map $M$ with a single $z$-gon we can define linear functions $i_P : \mathcal{E} \to \mathcal{E}$ and $\kappa_P : \mathcal{E} \to \mathcal{E}$ as follows. They are defined in the singletons and extended by linearity. Let $i_P(x)$ be the set of edges occurring once in the cyclic sequence $P$ between the two occurrences of $x$. Let $\kappa_P(x) = x$ if $x$ is traversed twice in the same direction in the zigzag path ($x$ is a black vertex in the medial map), and $\kappa_P(x) = \emptyset$, if $x$ is traversed in opposite direction in the zigzag path ($x$ is a white vertex in the medial map). Let $c_P = \kappa_P + i_P$. It is easy to verify that $c_P(x)$ is the set of edges occurring once in a closed path in $G_M$. Therefore, $c_P(x) \in \mathcal{V}^\perp$. In the above figure we see that $c_P(1) = \{1\} \cup \{2, 6, 7\}$, $c_P(7) = \emptyset \cup \{1, 4, 8, 9\}$. Indeed, $\{1, 2, 6, 7\}$ and $\{1, 4, 8, 9\}$ are members of $\mathcal{V}^\perp$. From the definitions, it follows that if $P$ has a single vertex, for any $x$, $\kappa_P(x) + \kappa_{P^\sim}(x) = x$ and that $c_{P^\sim}(x) + c_P(x) = x$ and $c_{P^\sim} + c_P$ is the identity linear transformation.

Let $x$ be a loop in an arbitrary map $M^t$. The loop is *balanced* if going around the $v$-gon corresponding to the vertex to which the loop is attached, and orienting the short edges accordingly the rectangle corresponding to $x$ gets short edges pointing in opposite directions. Otherwise the loop is *unbalanced*. The set of balanced loops of $M^t$ is denoted by $bal(M)$ and the set of its unbalanced loops is denoted by $unbal(M)$. Note that following the single zigzag in $G_M$ an edge $x$ is traversed twice in the same direction if and only if $x$ is a balanced loop in map $P^t$. In Fig. 4 we depict the situation in the rectangle $x$ of $M$ and $P$ corresponding to the edge $x$.

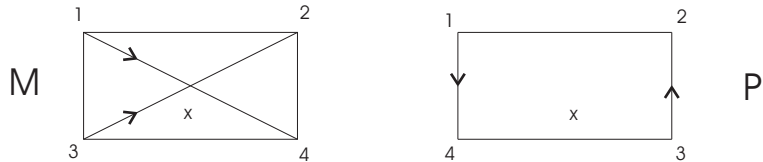

Fig. 4: The two passages through an edge following the zigzag viewed in the phial map

This observation shows that if $M$ is any map with a single $v$-gon, then $c_M(x) = i_M(x) + \kappa_M(x)$ is a cycle in $G_P$. Moreover, $\kappa_M(x) = x$ if $x$ is a balanced loop, otherwise $\kappa_M(x) = \emptyset$. The basic result about maps with a single zigzag is the following Theorem. The idea for its proof is taken from [Read and Rosensthiel, 1978] in a different context.

**Theorem 2** *If $M$ is a map with a single zigzag, $D$ its dual and $P$ its phial, then*

*(a)* $\mathrm{Im}(c_P) = \mathcal{V}^\perp$, *(b)* $\mathrm{Ker}(c_P) = \mathcal{V}$, *(c)* $\mathrm{Im}(c_{P^\sim}) = \mathcal{F}^\perp$, *(d)* $\mathrm{Ker}(c_{P^\sim}) = \mathcal{F}$.

**Proof:** From the signed sequence $P$ we to obtain a sequence of digraphs as follows. Start by drawing an oriented edge corresponding to first element in $P$. We proceed drawing oriented edges in the $P$-order, without lifting the pencil from the paper, as long as this is possible. Each time that an edge occurs for the first time we must draw it as a pendant edge (the final vertex must be a new one). The rule not to lift the pencil is not possible to obey when an edge $x$ occurring for the second time is not incident to the last vertex reached. At each such occurrence, we make a copy of the graph drawn so far, denoting it by $L_x$. Next, we make the necessary identification of two vertices in the new copy and proceed. Assuming $k$ identifications are necessary, we have defined a sequence of $k + 1$ graphs $L_{x_1}, L_{x_2}, \ldots, L_{x_k}, L_{x_{k+1}} = L$, where edge $x_i$ forced the $i$-th identification. If no identifications were made ($k = 0$), then we would have $m + 1$ vertices, where $m$ is the number of edges of $G_M$, which is the same as the number of edges of $L$. Each identification reduces the number of vertices by one. Therefore, if $v$ is the number of vertices of $L$ we have $m + 1 - k = v$ or $k = m + 1 - v$. Observe that this value of $k$ is the dimension of $CS(L)$, since $L$ is connected. We claim that $L = G_M$. Since the $P$-order can be realized by traversing edges in graph $G_M$, from the construction of $L$ as the graph with a minimum number of vertices where this is possible, it follows that $G_M$ can be obtained from $L$ by further identifications of vertices. However, every consecutive pair of elements in $P$ defines two edges which are consecutive in the zigzag of $M^t$. Therefore, they occur at a vertex of $G_M$ and no further identifications are necessary. It follows that $L = G_M$. If the two traversals of the edge $x_i$ are in opposite direction, the sequence of edges between the two occurrences of $x_i$ defines a reentrant path in $L_{i+1}$ which is not reentrant in $L_i$. The same is true for the set of edges between the two occurrences of $x_i$ plus edge $x_i$, in the case that $x_i$ is traversed in the same direction. From these facts we obtain that $c_P(x_i)$ is a cycle in $L_{i+1}$ and not a cycle in $L_i$. The set $\{c_P(x_1), \ldots, c_P(x_k)\}$ is linearly independent: indeed, a non-empty null linear combination of these vectors implies that the highest indexed one $c_P(x_h)$ is a sum of others with smaller indices. This is a contradiction because a $c_P(x_j)$ with $j < h$ is a cycle in $L_{x_h}$. Their sum would be a cycle in this graph, conflicting with the fact that $c_P(x_h)$ is not a cycle in it. Since $k$ is the dimension of the cycle space $\mathcal{V}^\perp$, $\{c_P(x_1), \ldots, c_P(x_k)\}$ is a basis for it. The proof of $(a)$ is complete.

To prove $(b)$, recall that over any field, $GF(2)$ in particular, $\mathcal{V} = BS(G_M) = (CS(G_M))^\perp$. In face of $(a)$, a subset of edges $A$ satisfies $A \in \mathcal{V}$ if and only if $|A \cap c_P(x)|$ is even for all edges $x$ of $G_M$. For $A, B \subseteq E(G_M)$, define the bilinear form $\langle A, B \rangle$ on $GF(2)$ by 0 if $|A \cap B|$ is even and 1 otherwise. We have $\langle A, c_P(x) \rangle = \sum_{a \in A} \langle a, c_P(x) \rangle = \sum_{a \in A} \langle c_P(a), x \rangle = \langle c_P(A), x \rangle$. It follows that $A \in \mathcal{V} \Leftrightarrow \langle c_P(A), x \rangle = 0$ for every edge $x$. The last condition is satisfied if and only if $c_P(A) = \emptyset$. This establishes $(b)$. To prove $(c)$ and $(d)$ apply $(a)$ and $(b)$ to map $D$, whose phial is $P^\sim$. ∎

The *connectivity of a closed surface* $S$ is denoted by $\xi(S)$ and is defined as $2 - \chi(S)$, where $\chi(S)$ is the Euler characteristic of $S$. The next result relates sum and intersection of bond and cycle spaces to the composition $c_{P^\sim} \circ c_P$. The dimension of these spaces are related to the connectivity of $\mathrm{Surf}(M, D)$.

**Theorem 3** *If $M$ is a map with one $z$-gon, then*

(a) $\dim[Im(c_{P\sim} \circ c_P)] = \xi(\text{Surf}(M,D))$,
(b) $Im(c_{P\sim} \circ c_P) = \mathcal{V}^\perp \cap \mathcal{F}^\perp$, (c) $\text{Ker}(c_{P\sim} \circ c_P) = \mathcal{V} + \mathcal{F}$.

**Proof:** By Theorem 1(a) we have $Im(c_P) = \mathcal{V}^\perp$. Therefore, $Im(c_{P\sim} \circ c_P) = Im(c_{P\sim}|\mathcal{V}^\perp)$, where, "$|\mathcal{V}^\perp$" stands for restriction to $\mathcal{V}^\perp$. By the fundamental theorem for homomorphisms [Godement, 1968] applied to $c_{P\sim}|\mathcal{V}^\perp$ we have

$$\dim[Im(c_{P\sim}|\mathcal{V}^\perp)] + \dim(\text{Ker}(c_{P\sim}|\mathcal{V}^\perp)] = \dim \mathcal{V}^\perp.$$

By part (b) of previous theorem, applied to map $P^\sim$, it follows that $\text{Ker}(c_P) = \mathcal{F}$. Since $\mathcal{F} \subseteq \mathcal{V}^\perp$, we get $\text{Ker}(c_{P\sim}|\mathcal{V}^\perp) = \mathcal{F}$. Therefore, $\dim[\text{Ker}(c_{P\sim}|\mathcal{V}^\perp)] = f - 1$. By the above equations, it follows that $\dim[Im(c_{P\sim} \circ c_P)] = (|E| - v + 1) - (f-1) = (|E| + 2) - (v+f) = \xi(\text{Surf}(M,D))$ and the proof of part (a) is complete.

To prove (b) we claim that $\dim(\mathcal{V}^\perp \cap \mathcal{F}^\perp) - \dim(\mathcal{V} \cap \mathcal{F}) = \xi(\text{Surf}(M,D))$. Note that $\dim(\mathcal{V}^\perp \cap \mathcal{F}^\perp) = (|E| - v + 1) + (|E| - f + 1) - \dim(\mathcal{V}^\perp + \mathcal{F}^\perp)$. Also that $\dim(\mathcal{V}^\perp + \mathcal{F}^\perp) + \dim(\mathcal{V} \cap \mathcal{F}) = |E|$, because they are orthogonal subspaces. Whence, $\dim(\mathcal{V}^\perp \cap \mathcal{F}^\perp) = (|E| + 2) - (v+f) + \dim(\mathcal{V} \cap \mathcal{F})$. But $(|E| + 2) - (v+f) = \xi(\text{Surf}(M,D))$, establishing the claim. By the Absorption Property, Theorem 1, it follows that $\mathcal{V} \cap \mathcal{F} = \mathcal{V} \cap \mathcal{F} \cap \mathcal{Z} = \{\emptyset\}$ (since $|VG_P| = 1$, $\mathcal{Z}$ is the null space). Thus, $\dim(\mathcal{V}^\perp \cap \mathcal{F}^\perp) = \xi(\text{Surf}(M,D))$. By part (a) $\dim[Im(c_{P\sim} \circ c_P)] = \xi(\text{Surf}(M,D))$. Note that $Im(c_P) = \mathcal{V}^\perp$ and $Im(c_{P\sim}) = \mathcal{F}^\perp$ imply the inclusion $Im(c_{P\perp} \circ c_P) \subseteq \mathcal{V}^\perp \cap \mathcal{F}^\perp$. Thus, $Im(c_P \circ c_{P\sim})$ is a subspace of $\mathcal{V}^\perp \cap \mathcal{F}^\perp$ which has the same dimension as the whole space. So, $Im(c_{P\sim} \circ c_P) = \mathcal{V}^\perp \cap \mathcal{F}^\perp$, concluding (b).

To prove (c) take $X \in \mathcal{V} + \mathcal{F}$; say that $X = U + W$ with $U \in \mathcal{V}$ and $W \in \mathcal{F}$. We have $(c_{P\sim} \circ c_P)(U + W) = c_{P\sim}(c_P(U) + c_P(W)) = c_{P\sim}(c_P(W))$. The latter equality follows because $c_P(U) = \emptyset$ by part (b) of previous Theorem. Since $c_P + c_{P\sim} = id$, they commute and we have $c_{P\sim} \circ c_P(W) = c_P \circ c_{P\sim}(W) = c_P(\emptyset) = \emptyset$. That $c_{P\sim}(W) = \emptyset$ follows from part (b) of Theorem 2 applied to map $P^\sim$. We conclude then that $\mathcal{V} + \mathcal{F} \subseteq \text{Ker}(c_{P\sim} \circ c_P)$. The value of $\dim[\text{Ker}(c_{P\sim} \circ c_P)]$ is $|E| - \xi(\text{Surf}(M,D))$. This follows from the fundamental theorem for homomorphisms and from part (a). Note that $\dim(\mathcal{V} + \mathcal{F}) = (v-1) + (f-1) = |E| - \xi(\text{Surf}(M,D))$. Since $\mathcal{V} + \mathcal{F}$ is a subspace of $\text{Ker}(c_{P\sim} \circ c_P)$ and has the same dimension, it follows that it is equal to $\text{Ker}(c_{P\sim} \circ c_P)$. This concludes (c). ∎

## 4 Maps with a single zigzag and a single face

Finally, we prove a result on maps $M^t$ having a single zigzag **and** a single face. By subdividing some edges of an arbitrary graph $G$ it is possible to embed $G$ in some surface so that the resulting $t$-map has a single face and a single zigzag. Below we present an embedding of $K_4$ having this property.

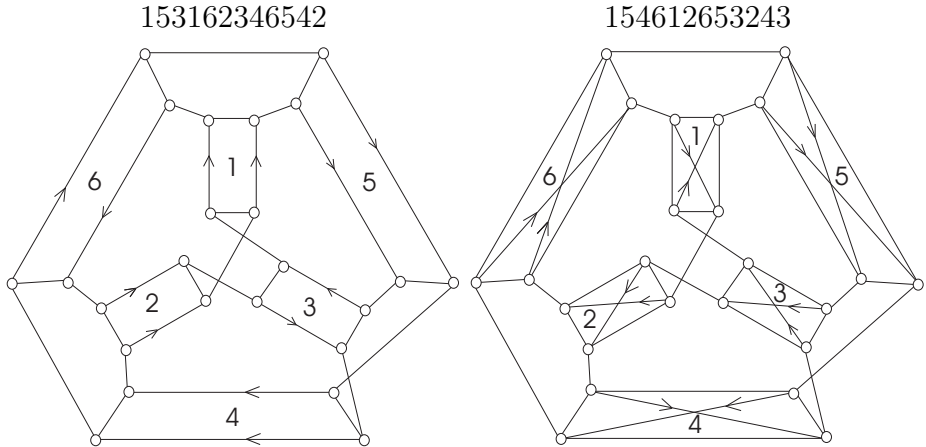

Fig. 5: An embedding of $K_4$ with a single face and a single zigzag

We have the following values for the functions $c_{P\sim}$ and $c_D$ on the singletons:

$c_{P\sim}(1) = \{1,5,4,6\}$ $c_D(1) = \{5,3\}$ $c_{P\sim}(2) = \{6,5,3\}$ $c_D(2) = \{3,6,5\}$

$c_{P\sim}(3) = \{2,4\}$ $c_D(3) = \{3,1,6,2\}$ $c_{P\sim}(4) = \{4,3,1,5\}$ $c_D(4) = \{6,5\}$

$c_{P\sim}(5) = \{4,1,2\}$ $c_D(5) = \{4,2,1\}$ $c_{P\sim}(6) = \{1,2\}$ $c_D(6) = \{6,2,3,4\}$

Their composition $c_{P\sim} \circ c_D$ is the identity:

$c_{P\sim} \circ c_D(1) = c_{P\sim}(\{5,3\}) = \{4,1,2\} + \{2,4\} = \{1\}$

$c_{P\sim} \circ c_D(2) = c_{P\sim}(\{3,6,5\}) = \{2,4\} + \{1,2\} + \{4,1,2\} = \{2\}$

$c_{P\sim} \circ c_D(3) = c_{P\sim}(\{3,1,6,2\}) = \{2,4\} + \{1,5,4,6\} + \{1,2\} + \{6,5,3\} = \{3\}$

$c_{P\sim} \circ c_D(4) = c_{P\sim}(\{6,5\}) = \{4,1,2\} + \{1,2\} = \{4\}$

$c_{P\sim} \circ c_D(5) = c_{P\sim}(\{4,2,1\}) = \{4,3,1,5\} + \{6,5,3\} + \{1,5,4,6\} = \{5\}$

$c_{P\sim} \circ c_D(6) = c_{P\sim}(\{6,2,3,4\}) = \{1,2\} + \{6,5,3\} + \{2,4,\} + \{4,3,1,5\} = \{6\}$

As our final Theorem show, this is a general property of $t$-maps having a single face and a single zigzag.

**Theorem 4** *If $M^t$ has a single face and a single zigzag, then $c_{P\sim} \circ c_D$ is the identity on $\mathcal{E}$.*

**Proof:** We prove that $c_{P\sim} \circ c_D(x) = x$ in the three following cases. $I : x \in bal(D); II : x \in unbal(D), x \in bal(P^\sim); III : x \in unbal(D), x \in unbal(P^\sim)$.

Assuming the hypothesis of case $I$, consider the map $D' = D(s\ell, x)$ and its phial $P' = P^\sim(\ell d, x)$. The unique vertex of $G_D$ breaks into two in $G_{D'}$, arising the bond $x \cup i_D(x)$ linking these two vertices. Therefore, by Theorem 1(d), $c_{P'}(x \cup i_D(x)) = \emptyset$. Note that $c_{P'}(x) = x + c_{P^\sim}(x)$ and $c_{P'}(y) = c_{P^\sim}(y)$, for $y \neq x$. So, we get $x + c_{P^\sim}(x) + c_{P^\sim}(i_D(x)) = \emptyset$. But $i_D(x) = x + c_D(x)$, since $x \in bal(D)$. So, $x + c_{P^\sim}(x) + c_{P^\sim}(x) + c_{P^\sim}(c_D(x)) = \emptyset$, or $c_{P^\sim}(c_D(x)) = x$, establishing case $I$.

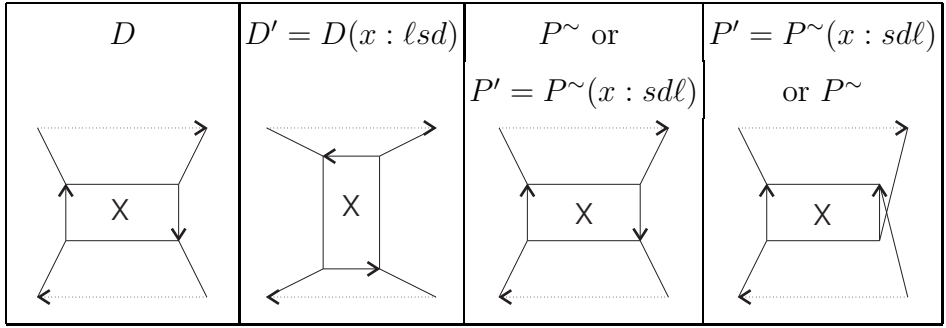

Fig. 6: Maps involved in the proof of Case I.

Assume the hypotheses of case $II$, and consider map $D'' = D(sd, x)$, and its phial $P'' = P^\sim(sd; x)$. The set $x \cup i_D(x)$ is the unique bond of $G_{D''}$, linking its two vertices. Therefore, by Theorem 1(d), $c_{P''}(x \cup i_D(x)) = \emptyset$. Observe since $x \in bal(P'')$, $c_{P''}(x) = c_{P^\sim}(x) = x \cup i_{P^\sim}(x)$. For $y \notin i_{P^\sim}(x)$, $i_{P''}(y) = i_{P^\sim}(y)$ and $\kappa_{P''}(y) = \kappa_{P^\sim}(y)$. For $y \in i_{P^\sim}(x)$, $i_{P''}(y) = y + i_{P^\sim}(y) + i_{P^\sim}(x)$ and $\kappa_{P''}(y) = y + \kappa_{P^\sim}(y)$. We get $\emptyset = c_{P''}(x + i_D(x)) = c_{P''}(x) + c_{P''}(i_D(x)) = c_{P^\sim}(x) + \sum\{c_{P''}(y) \mid y \in i_D(x) \setminus i_{P^\sim}(x)\} + \sum\{c_{P''}(y) \mid y \in i_D(x) \cap i_{P^\sim}(x)\} = c_{P^\sim}(x) + \sum\{c_{P''}(y) \mid y \in i_D(x) \setminus i_{P^\sim}(x)\} + \sum\{i_{P''}(y) + \kappa_{P''}(y) \mid y \in i_D(x) \cap i_{P^\sim}(x)\} = c_{P^\sim}(x) + \sum\{c_{P^\sim}(y) \mid y \in i_D(x) \setminus i_{P^\sim}(x)\} + \sum\{c_{P^\sim}(y) + i_{P^\sim}(x) \mid y \in i_D(x) \cap i_{P^\sim}(x)\} = c_{P^\sim}(x) + c_{P^\sim}(i_D(x)) + m i_{P^\sim}(x)$, where $m = |i_D(x) \cap i_{P^\sim}(x)|$. Assuming that $m$ is odd, we have $c_{P^\sim}(x) + c_{P^\sim}(i_D(x)) = i_{P^\sim}(x) = c_{P^\sim}(x) + x$. Therefore, $c_{P^\sim}(c_D(x)) = c_{P^\sim}(i_D(x)) = x$, establishing case $II$, provided $m$ is odd. If it is even we have proved that $c_{P^\sim}(x) + c_{P^\sim}(i_D(x)) = \emptyset$. This implies that $x + i_D(x)$ is in the kernel of $c_{P^\sim}$, which is $\mathcal{F}^\perp = \{\emptyset\}$. It follows that $x + i_D(x) = \emptyset$, or $x = i_D(x)$, a contradiction because $x \notin i_D(x)$. The proof of case $II$ is complete.

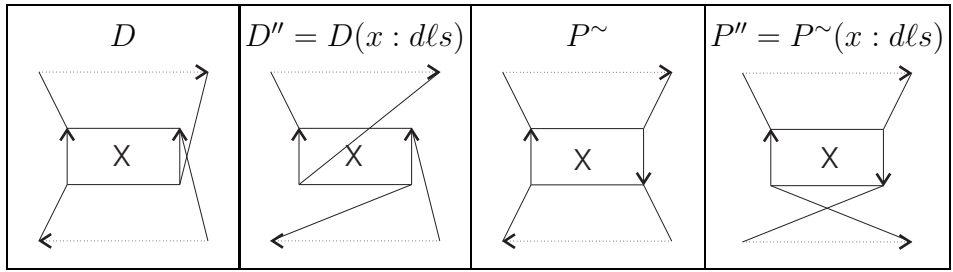

Fig. 7: Maps involved in the proof of Case II.

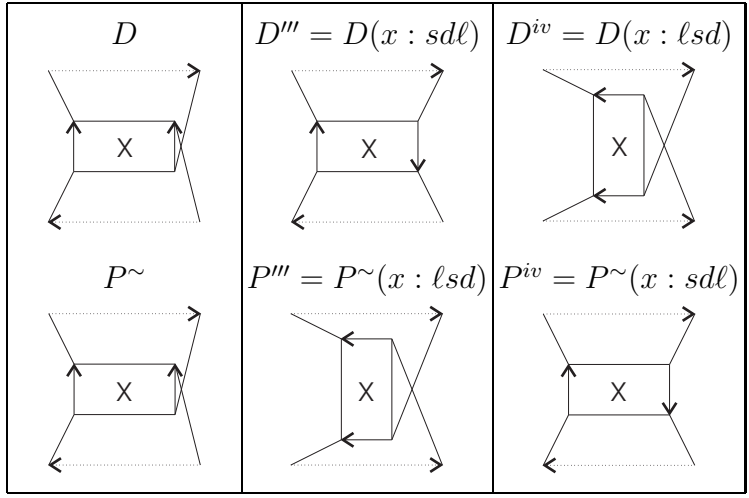

Fig. 8 : Maps involved in the proof of Case III.

Assume the hypotheses of case $III$, let $D''' = D(\ell d, x)$ and $P''' = P^\sim(s\ell, x)$ be its phial. Both $D'''$ and $P'''$ have a single vertex and $x \in bal(D''')$. Therefore, by case $I$, $c_{P'''}(c_{D'''}(x)) = x$. As $c_{D'''}(x) = x + c_D(x) = x + i_D(x)$, we get the equality (E) $c_{P'''}(x) + c_{P'''}(i_D(x)) = x$. Since $x \in unbal(P''')$, $c_{P'''}(x) = c_{P^\sim}(x) = i_{P^\sim}(x)$. Expressing $c_{P'''}(i_D(x))$ in terms of $c_{P^\sim}$: $c_{P'''}(i_D(x)) = c_{P'''}(i_D(x) \setminus i_{P^\sim}(x)) + c_{P'''}(i_D(x) \cap i_{P^\sim}(x))$. As for $P''$ in case $II$, we get $c_{P'''}(i_D(x) \setminus i_{P^\sim}(x)) = c_{P^\sim}(i_D(x) \setminus i_{P^\sim}(x))$ and $c_{P'''}(i_D(x) \cap i_{P^\sim}(x)) = c_{P^\sim}(i_D(x) \cap i_{P^\sim}(x)) + |i_D(x) \cap i_{P^\sim}(x)| i_{P^\sim}(x)$. Whence, $c_{P'''}(i_D(x)) = c_{P^\sim}(i_D(x)) + m i_{P^\sim}(x)$, where $m = |i_D(x) \cap i_{P^\sim}(x)|$. Rewriting (E) we get $c_{P^\sim}(x) + c_{P^\sim}(i_D(x)) + m i_{P^\sim}(x) = x$. As we show shortly, $m$ is odd, and so, $c_{P^\sim}(x) + c_{P^\sim}(i_D(x)) + i_{P^\sim}(x) = x$. But in case $III$ $c_{P^\sim}(x) = i_{P^\sim}(x)$ and $i_D(x) = c_D(x)$. It follows that $c_{P^\sim}(c_D(x)) = x$, establishing this final case, provided $m$ is odd. Consider the map $D^{iv} = D(s\ell, x)$ and its phial $P^{iv} = P^\sim(\ell d, x)$. They both have a single vertex and satisfy the hypotheses of case $II$, and so we have proved that $|i_{D^{iv}}(x) \cap i_{P^{iv}}(x)|$ is odd. To conclude the proof, just note that $i_{D^{iv}}(x) = i_D(x)$ and $i_{P^{iv}}(x) = i_{P^\sim}(x)$. Thus, $m = |i_D(x) \cap i_{P^\sim}(x)| = |i_{D^{iv}}(x) \cap i_{P^{iv}}(x)|$ is odd. ∎